\documentclass[11pt]{article}
\renewcommand{\baselinestretch}{1.08}\rm
\newcommand{\conv}{\mbox{\rm conv}}

\newcommand{\inter}{\mbox{\rm int}}
\newcommand{\relinter}{\mbox{\rm rel-int}}
\newcommand{\vol}{\mbox{\rm vol}}
\newcommand{\ca}{\mbox{\rm ca}}
\author{Oleg Pikhurko\thanks{Supported by a Senior Rouse Ball
 Studentship, Trinity College, Cambridge, UK.}\\
 DPMMS, Centre for Mathematical Sciences\\
 Cambridge University, Cambridge~CB3~0WB, England\\
 E-mail: {\tt O.Pikhurko@dpmms.cam.ac.uk}}

\usepackage{amssymb}
\newcommand{\eqref}[1]{\mbox{\rm(\ref{#1})}}
\newcommand{\pl}{}

\newcommand{\B}[1]{{\bf #1}}
\newcommand{\I}[1]{{\mathbb #1}}

\newcommand{\comment}[1]{}
\newcommand{\qed}{\nolinebreak\mbox{\hspace{5 true pt}%
\rule[-0.85 true pt]{3.9 true pt}{8.1 true pt}}}

\newtheorem{theorem}{Theorem}
\newtheorem{lemma}[theorem]{Lemma}

\begin{document}

\renewcommand{\equiv}{:=}

\title{Lattice Points in Lattice Polytopes}

\author{Oleg Pikhurko\thanks{This research was carried out during the
author's stay at the Technical University, Berlin, sponsored by the
German Academic Exchange Service (DAAD) and the Rouse Ball Travelling Fund of
Trinity College, Cambridge.}\\
 DPMMS, Centre for Mathematical Sciences\\
 Cambridge University, Cambridge~CB3~0WB, England\\
 E-mail: {\tt O.Pikhurko@dpmms.cam.ac.uk}}


\maketitle

\newcommand{\pcadl}{8d\cdot (8l+7)^{2^{2d+1}}}
\newcommand{\scadl}{8\cdot (8l+7)^{2^{d+1}}}
\newcommand{\smdl}{(8l+7)^{-2^{d+1}}/8}
\newcommand{\sdkl}{2^{3d-2}\cdot l^d\cdot (8l+7)^{(d-1)2^{d+1}}\cdot k/d!}
\newcommand{\pdkl}{(8dl)^d\cdot (8l+7)^{d\cdot 2^{2d+1}}\cdot k}

\begin{abstract}
 We show that, for any lattice polytope $P\subset\I R^d$, the set
$\inter(P)\cap l\I Z^d$ (provided it is non-empty) contains a point whose
coefficient of asymmetry with respect to $P$ is at most $\pcadl$. If,
moreover, $P$ is a simplex, then this bound can be improved to $\scadl$.

As an application, we deduce new upper bounds on the volume of a
lattice polytope, given its dimension and the number of sublattice points
in its interior.
\end{abstract}

\section{Introduction}\label{intro}

A {\em lattice polytope} in $\I R^d$ is a convex polytope whose
vertices are lattice points, that is, points in $\I Z^d$. For an
integer $l\ge 1$, let $I_l(P)=\inter(P)\cap l\I Z^d$ be the set of
interior points of $P$ whose coordinates are integers divisible by
$l$.

Of course, some points of $I_l(P)$ can lie `close' to $\partial P$,
the boundary of $P$. However, our Theorem~\ref{th:\pl:ca-polytope} shows
that, provided $I_l(P)\not=\emptyset$, there is $\B w\in I_l(P)$
with
 \begin{equation}\label{eq:\pl:ca-polytope-i}  \ca(\B w,P)\le\pcadl,\end{equation}
 where $\ca(\B w,P)$ is the {\em coefficient of asymmetry} of $P$ about
$\B w$:
 $$\ca(\B w,P)=\max_{|\B y|=1} \frac{\max\{\lambda\mid \B w+
\lambda \B y\in P\}}{\max\{\lambda\mid \B w- \lambda \B y\in P\}}.$$
Although the function in the right-hand side of~\eqref{eq:\pl:ca-polytope-i} is
enormous, the main point is that it depends only on $d$ and $l$.
 
We prove an inequality of this type for the case of a simplex
$S$ first. Namely, Theorem~\ref{th:\pl:ca-simplex} implies that, for some
$\B w\in I_l(S)$, 
 \begin{equation}\label{eq:\pl:ca-simplex-i} \ca(\B w,S)\le\scadl.\end{equation}
 Here the claim essentially
concerns the barycentric coordinates $(\alpha_0,\dots,\alpha_d)$ of
$\B w$ inside $S$ because of the easy relation
 \begin{equation}\label{eq:\pl:ca-m}\ca(\B w,S)=\max_{0\le i\le d}\frac{1-\alpha_i}{\alpha_i}
=\frac{1}{m_S(\B w)}-1,\end{equation}
 where $m_S(\B w)\equiv \min_{0\le i\le d} \alpha_i$ is
the smallest barycentric coordinate of $\B w\in S$. Define
 \begin{equation}\label{eq:\pl:beta} \beta(d,l)\equiv\inf_S \max\{m_S(\B w)\mid \B w\in I_l(S)\},\end{equation}
 where the infimum is taken over all lattice simplices $S$ with
$I_l(S)\not=\emptyset$. (For example, it is easy to see that
$\beta(1,l)=\frac1{l+1}$.) Thus we have to prove a positive lower bound
on $\beta(d,l)$. The gist the proof is
that if we have  $\B w\in I_l(S)$ with $m_S(\B
w)$ being `small', then using one approximation lemma of Lagarias and
Ziegler~\cite{lagarias+ziegler:91} we can `jump' to another vertex $\B
w'\in I_l(S)$ with $m_S(\B w')>m_S(\B w)$, see
Theorem~\ref{th:\pl:ca-simplex}.

In fact, one result of Lawrence~\cite[Theorem~3]{lawrence:91} implies
that $\beta(d,1)>0$ but does not give any explicit bound, see
Section~\ref{lawrence} here.

It would be interesting to know how far our bounds~\eqref{eq:\pl:ca-polytope-i}
and~\eqref{eq:\pl:ca-simplex-i} are from the best possible values. The best
values that we know arise from the following family of lattice simplices.

Define inductively the sequence $t_{d,l}$ by $t_{1,l}=l+1$ and
$t_{d+1,l}=t_{d,l}^2 -t_{d,l}+1$. (This sequence appears
in~\cite{lagarias+ziegler:91}.) Consider the simplex
 \begin{equation}\label{eq:\pl:Bdkl} B_{d,l}\equiv\conv\{t_{1,l}\,\B e_1,\,\dots,\,
t_{d-1,l}\,\B e_{d-1},\, t_{d,l}\,\B e_d,\, -\B e_d\}\subset \I
R^d,\end{equation}
 where $(\B e_1,\dots,\B e_d)$ is the standard basis.  It is not hard
to see that $I_l(B_{d,l})=\{l\B 1, l(\B 1-\B e_d)\}$,
cf.~\cite[Proposition~2.6]{lagarias+ziegler:91}.  We have
$m_{B_{d,l}}(l\B 1)=m_{B_{d,l}}(l(\B 1-\B e_d))=\frac l{t_{d,l}^2-1}$:
the vertex $l\B 1=(l,\dots,l)$, for example, has barycentric
coordinates
 $$\left(\frac l{t_{1,l}},\ \dots,\ \frac l{t_{d-1,l}},\ \frac
l{t_{d,l}-1}\times \frac{t_{d,l}}{t_{d,l}+1},\ \frac
l{t_{d,l}-1}\times \frac{1}{t_{d,l}+1}\right).$$
 One can show that $t_{d,l}\ge (l+1)^{2^{d-2}}+1$ for $d\ge 2$ by
considering $u_{d,l}=t_{d,l}-1$; hence
 \begin{equation}\label{eq:\pl:beta-upper} \beta(d,l)\le \frac{l}{t_{d,l}^2-1}\le
l(l+1)^{- 2^{d-1}},\quad d\ge 2.\end{equation} 

Thus~\eqref{eq:\pl:ca-simplex-i} establishes the correct type of dependence on
$d$ and $l$, although the gap between the bounds is huge.  Perhaps,
$B_{d,l}$ gives the actual value of the function $\beta(d,l)$ as well
as the sharp bound for~\eqref{eq:\pl:ca-polytope-i}.

To extend Theorem~\ref{th:\pl:ca-simplex} to a general lattice polytope
$P\subset \I R^d$ (Theorem~\ref{th:\pl:ca-polytope}), we try to find a lattice
polytope $P'\subset P$ with few vertices such that
$I_l(P')\not=\emptyset$ and a homothetic copy of $P'$ covers $P$.  The
latter condition gives an upper bound on $\ca(\B w,P)$ in terms of
$\ca(\B w,P')$ for $\B w\in\inter(P')$, see Lemma~\ref{lm:\pl:ca}, and is
satisfied if, for example, $P'\supset S$, where $S\subset P$ is a
simplex of the maximum volume.  But to get a non-empty $I_l(P')$ we
may have to add as many as $d$ extra vertices to $S$. It is now
possible to define our jumps within $P'$ to get the required $\B w\in
I_l(P')$. However, the bound~\eqref{eq:\pl:ca-polytope-i} for $d$-polytopes
that we obtain is comparable with that for $2d$-dimensional simplices;
we believe that we lose here too much but we have not found any
better argument.

Next, we investigate the following problem. Let $p(d,k,l)$ (resp.\
$s(d,k,l)$) be the maximum volume of a lattice polytope (resp.\
simplex) $P\subset\I R^d$ with $|I_l(P)|=k$. As for any $d\ge2$ there
exist lattice simplices with no lattice points in the interior and of
arbitrarily large volume, we restrict our consideration to the case
$k\ge1$.

Trivially, $p(1,k,l)=s(1,k,l)=(k+1)l$. A result of
Scott~\cite{scott:76} implies that $p(2,1,1)=s(2,1,1)=\frac92$ and
$p(2,k,1)=s(2,k,1)=2(k+1)$ for $k\ge
2$. Hensley~\cite[Theorem~3.6]{hensley:83} showed that $p(d,k,1)$
exists (i.e., it is finite) for $k\ge 1$. The method of Hensley was
sharpened by Lagarias and
Ziegler~\cite[Theorem~1]{lagarias+ziegler:91}, who showed that
 \begin{equation}\label{eq:\pl:LZ91} p(d,k,l)\le kl^d (7(kl+1))^{d2^{d+1}},\end{equation}
 and also observed that, for any fixed $(d,k,l)$, there are finitely
many (up to a $GL_n(\I Z)$-equivalence) lattice polytopes $P\subset \I
R^d$ with $|I_l(P)|=k\ge 1$.

Lagarias and Ziegler~\cite[Theorem~2.5]{lagarias+ziegler:91} proved
the following extension 
of a theorem of Mahler~\cite{mahler:39}: ``A convex body $K\subset\I
R^d$ with $k=|I_l(K)|\ge 1$ satisfies
 \begin{equation}\label{eq:\pl:M39} \vol(K)\le \left(l(\ca(\B w,K)+1)\right)^d\cdot k,\end{equation}
 for any $\B w\in I_l(K)$.''

Combining~\eqref{eq:\pl:M39} with~\eqref{eq:\pl:ca-polytope-i} (or more exactly
with~\eqref{eq:\pl:ca-polytope-th}), we obtain that
 \begin{equation}\label{eq:\pl:p}p(d,k,l)\le\pdkl.\end{equation}

A theorem of Blichfeldt~\cite{blichfeldt:14} says that $|P\cap\I
Z^d|\le n+n!\,\vol(P)$; combined with~\eqref{eq:\pl:p} it gives an
upper bound on $|P\cap\I Z^d|$ in terms of $|I_l(P)|$ (if the latter
set is non-empty).

An upper bound on $s(d,k,l)$ can be obtained by applying~\eqref{eq:\pl:M39}
to~\eqref{eq:\pl:ca-simplex-i}. However, we obtain a better bound in
Theorem~\ref{th:\pl:sdkl} by exploiting the geometry of a simplex, namely we
show that
 \begin{equation}\label{eq:\pl:sdkl-i} s(d,k,l)\le\sdkl.\end{equation}

The best lower bound on $p(d,k,l)$ and $s(d,k,l)$ that we know (except
for $(d,k,l)=(2,1,1)$), comes from the consideration of the simplex
 $$S_{d,k,l}\equiv\conv\{\B 0,\,t_{1,l}\, \B e_1,\,\,\dots,\,t_{d-1,l}\,\,\B
e_{d-1},(k+1)(t_{d,l}-1)\B e_d\},$$
 which satisfies $I_l(S_{d,k,l})=\{l\B 1+il\B e_d\mid 0\le i\le
k-1\}$, see~\cite[Proposition~5.6]{lagarias+ziegler:91}. This
demonstrates that
 $$s(d,k,l)\ge \vol(S_{d,k,l})> \frac{k+1}{d!\cdot l}\,(l+1)^{2^{d-1}},$$
 see formula~(2.13) in~\cite{lagarias+ziegler:91}. The family
$(S_{d,k,1})$ was found by Zaks, Perles and
Wills~\cite{zaks+perles+wills:82} and its generalization (the addition
of parameter $l$)---by Lagarias and
Ziegler~\cite{lagarias+ziegler:91}.

Again, we have the correct type of dependence of $d$, $k$ and $l$ but
the gap between the known bounds is huge. The ultimate aim would be to
find exact values, which is probably not hopeless because the above
contructions, believed to be extremal, are rather simple.

\section{Jumping inside a simplex}

We will use the following lemma of Lagarias and
Ziegler~\cite[Lemma~2.1]{lagarias+ziegler:91}.

\begin{lemma}\label{lm:\pl:LZ91} For a real $\lambda\ge 1$ and integer $n\ge 1$, define\begin{equation}\label{eq:\pl:LZ4}
 \delta(n,\lambda)=\left(7(\lambda+1)\right)^{-2^{n+1}}.\end{equation}
 Then, for all positive real
numbers $\alpha_1,\dots,\alpha_n$ satisfying
 $$1-\delta(n,\lambda)<\sum_{i=1}^n \alpha_i\le 1,$$
 there exist non-negative integers $P_1,\dots,P_n,Q$ such that
 \begin{eqnarray}
 Q&=& P_1+\cdots+P_n\ >\ 0,\label{eq:\pl:LZ1}\\
 \alpha_i&>&\frac{\lambda P_i}{\lambda Q+1} \mbox{ \ for $1\le i\le n$},\label{eq:\pl:LZ2}\\
 \lambda Q+1&\le & \delta(n,\lambda)^{-1}.\qed\label{eq:\pl:LZ3}
\end{eqnarray}\end{lemma}

The above lemma is the main ingredient in our `jumps.' Here it is applied
with $\lambda=\frac87\, l$. There is nothing special about the constant
$\frac87$ except it makes~\eqref{eq:\pl:LZ4} look simpler; any fixed number greater
than $1$ would do as well.

\begin{theorem}\label{th:\pl:ca-simplex} Let $l\ge 1$ and let $S=\conv\{\B
v_0,\dots,\B v_d\}\subset\I R^n$ be a lattice simplex. If
$\relinter(S)\cap l\I Z^n$ is non-empty, then it contains a point $\B
w$ with
 \begin{equation}\label{eq:\pl:ca-simplex} m_S(\B w)\ge \gamma\equiv
\delta(d,{\textstyle\frac87}\,l)/8=\smdl.\end{equation}\end{theorem} 
 \smallskip{\it Proof.}  We may assume that $n=d$ because we can always find a linear
transformation preserving the lattice $\I Z^n$ (and so $l\I Z^n$ as well) and
mapping $S$ into $\I R^d\subset \I R^n$.

Let $\B w=\sum_{i=0}^d \alpha_i \B v_i\in I_l(S)$, $\sum_{i=0}^d
\alpha_i =1$, be a vertex maximizing $m_S(\B w)$. Suppose that the
claim is not true. Assume that $\alpha_0\le\cdots\le\alpha_d$; then
$m_S(\B w)=\alpha_0<\gamma$. Let $j$ be the index with $\alpha_j<
8\gamma\le \alpha_{j+1}$; note that $j\le d-1$ is well-defined as
$\alpha_d\ge\frac1{d+1} \ge 8\gamma$.

We have $\sum_{i=0}^j \alpha_j < 8\gamma(j+1)$ which, as it is easy to see, does
not exceed $\delta(d-j,\frac87\, l)$ for $j\in[0,d-1]$. Hence,
Lemma~\ref{lm:\pl:LZ91} is applicable to the $d-j$ numbers
$\alpha_{j+1},\dots,\alpha_d$ and yields integers $P_{j+1},\dots, P_d,Q$
satisfying~\eqref{eq:\pl:LZ1}--\eqref{eq:\pl:LZ3}.

Consider the vertex
 $$\B w'=(lQ+1)\B w-\sum_{i=j+1}^d l P_i\B v_i\ \in\ l\I Z^d.$$
 We have $\B w'=\sum_{i=0}^d \alpha_i' \B v_i$, where, for
$i\in[0,j]$, $\alpha_i'\equiv (lQ+1)\alpha_i>\alpha_0$ and, for
$i\in[j+1,d]$, $\alpha_i'\equiv (lQ+1)\alpha_i-lP_i > \alpha_i/8\ge
\alpha_0$ by~\eqref{eq:\pl:LZ2}. As
 $$\sum_{j=0}^n\alpha_j'=(lQ+1)\sum_{i=0}^d\alpha_i-\sum_{i=j+1}^d
l P_i=1,$$
 the lattice point $\B w'$ lies in the interior of $S$ and
contradicts the choice of~$\B w$.\qed \medskip

\smallskip\noindent{\bf Remark.}  For $n=d$ the inequality~\eqref{eq:\pl:ca-simplex-i} claimed in the
introduction follows by applying~\eqref{eq:\pl:ca-m} to the vertex $\B w\in
I_l(S)$ given by Theorem~\ref{th:\pl:ca-simplex}\smallskip

\section{$\beta(d,l)$ for small $d$ and $l$}\label{small-beta}

Let us try to deduce some estimates of $\beta(d,l)$ when $d$ and $l$
are small. We have a general upper bound~\eqref{eq:\pl:beta-upper} which in
particular says that
 \begin{eqnarray}
 \beta(2,l)&\le&\frac l{t_{2,l}^2-1}\ =\
\frac{1}{(l+1)(l^2+l+2)},\label{eq:\pl:b2l}\\ 
 \beta(3,1)&\le& \frac1{48},\quad \beta(3,2)\ \le\
\frac{1}{924}.\label{eq:\pl:b31}\end{eqnarray}

Here we present some results obtained with the help of computer showing
that~\eqref{eq:\pl:b2l} and~\eqref{eq:\pl:b31} are probably sharp.

How can one get a lower bound on, for example, $\beta(2,1)$? Our approach
is the following.

Given a lattice simplex $S=\conv\{\B v_0,\B v_1,\B v_2\}\subset \I
R^2$, let $\B w$ be a lattice vertex maximizing $m_S$ over
$I_1(S)\not=\emptyset$. Write the barycentric representation $\B
w=\sum_{i=0}^2 \alpha_i\B v_i$. We have
 \begin{equation}\label{eq:\pl:b21-1} \alpha_0+\alpha_1+\alpha_2=1.\end{equation}
 Without loss of generality we may assume that
 \begin{equation}\label{eq:\pl:b21-2} 0\le \alpha_0\le \alpha_1\le \alpha_2.\end{equation}
 
Consider the vertex $\B w'=2\B w-\B v_3$ which is the jump of $\B w$
corresponding to $\B P=(0,0,1)$. Its barycentric coordinates
$(\alpha_0', \alpha_1',\alpha_2')=(2\alpha_0,2\alpha_1,2\alpha_2-1)$
satisfy $\alpha_1'\ge \alpha_0' >\alpha_0$. By the choice of
$\B w$, we must have
 \begin{equation}\label{eq:\pl:b21-3} 2\alpha_2-1\le \alpha_0.\end{equation}
 
Similarly, the $(0,1,1)$-jump $\B w'=3\B w-\B v_1-\B v_2$ satisfies
$\alpha_0'=3\alpha_0>\alpha_0$ and, in the view of $\alpha_2'\ge
\alpha_1'$, we obtain
 \begin{equation}\label{eq:\pl:b21-4} 3\alpha_1-1\le \alpha_0.\end{equation}

Not everything goes so smoothly if we consider e.g.\ the
$(0,1,2)$-jump, when we can only deduce that
 \begin{equation}\label{eq:\pl:b21-5} 4\alpha_1-1\le\alpha_0 \quad\mbox{or}\quad 4\alpha_2-2\le
\alpha_0.\end{equation}

How small can $\alpha_0$ be, given only the
constraints~\eqref{eq:\pl:b21-1}--\eqref{eq:\pl:b21-5} (which can be realized as a mixed
integer program)? Solving this MIP, we obtain
that $\alpha_0\ge\frac2{19}$. Knowing this bound, we can enlarge our
arsenal of jumps. For example, the vertex $\B w'=12\B w -\B v_1-4\B
v_2- 6\B v_2$ satisfies $\alpha_0'=12\alpha_0-1>\alpha_0$; hence
 \begin{equation}\label{eq:\pl:b21-6} 12\alpha_1-4\le \alpha_0\quad \mbox{or} \quad
12\alpha_2-6 \le \alpha_0.\end{equation}

The addition of~\eqref{eq:\pl:b21-6} to the system~\eqref{eq:\pl:b21-1}--\eqref{eq:\pl:b21-5},
improves our lower bound to $\alpha_0\ge\frac2{17}$.

In this manner we can repeatedly add new constraints to our MIP as
long as this improves the lower bound on $\alpha_0$. This was realized
as a program in C which can be linked with either CPLEX (commercial)
or {\tt lp\_solve} (public domain) MIP solver. The source
code is freely available~\cite{pikhurko:00:co} and the reader is welcome
to experiment with it.

When one runs the program, it seems that the obtained lower bound $f$
on $a_0$ approaches some limit without attaining it. And, of course,
the more iterations we perform, the larger coefficients in the
added inequalities are and the problem becomes more complex.
Hence, the question to what extend we can trust the output should be
considered. The CPLEX has the mechanism to set up various tolerances
which specify how far CPLEX allows variables to violate the bounds and
to be still considered feasible during perturbations. The manual does
not specify the guaranteed accuracy. As the coefficient at
$\alpha_0$ is always $1$, we assume that the error does not exceed
$\Delta$, the largest sum of absolute values of coefficients in one
inequality multiplied by the tolerance, which we set to $10^{-9}$, the
smallest value allowed by CPLEX's manual.

In particular, we allow a $(P_0,P_1,P_2)$-jump when we can guarantee
that $(1+P_0+P_1+P_2)\alpha_0 -P_0\ge \alpha_0$, namely, when $(f -
\Delta) (P_0+P_1+P_2) \ge P_0$.

We ran the program, with CPLEX~6.6, for various $d$ and $l$;
Table~\ref{beta-cplex} records the obtained lower bounds up to
$10^{-6}$. Unfortunately, we had no success when $d\ge 4$ or when
$d=3$ and $l\ge 3$ or when $d=2$ and $l\ge 28$: the obtained lower
bound was still zero when the MIP became too large to solve.

\begin{table}[p]\caption{Computed lower bounds on $\beta(d,l)$ using CPLEX.}
\label{beta-cplex}

\begin{tabular}{|c|c|c|c|c|c|c|}
\hline\hline

$d$ & $l$ & our upper & number of & output of & estimated & `guaranteed'\\
& & bound & iterations & the program & error & lower bound \\
\hline\hline
2 & 1 & 0.125000 & 375 & 0.124906 & $< 2\times 10^{-6}$ & 0.124904 \\
\hline
2 & 2 & 0.041667 & 167  & 0.041648 & $< 4\times 10^{-6}$ & 0.041644\\
\hline
2 & 3 & 0.017858 & 72 & 0.017850 & $< 5\times 10^{-6}$ & 0.017845\\
\hline
2 & 4 & 0.009091 & 38  & 0.009086 & $< 4\times 10^{-6}$ & 0.009082\\
\hline
2 & 5 & 0.005209 & 33 & 0.005205 & $<4\times 10^{-6}$ & 0.005201\\
\hline
2 & 6 &0.003247 & 40 & 0.003245 & $<4\times 10^{-6}$ & 0.003241\\
\hline
2 & 7 & 0.002156 & 49 & 0.002150 & $<1\times 10^{-6}$ & 0.002149\\ 
\hline
2 & 8 & 0.001502 & 65 & 0.001499 & $< 2\times 10^{-6}$ & 0.001497\\
\hline
2 & 9 & 0.001087 & 82 & 0.001085 & $< 2\times 10^{-6}$&0.001083\\
\hline
2 & 10 & 0.000812 & 97 & 0.000810 & $< 2\times 10^{-6}$&0.000808\\
\hline
2 & 11 & 0.000622 & 116 & 0.000621 & $< 2\times 10^{-6}$&0.000619\\
\hline
2 & 12 & 0.000487 & 144 & 0.000486 & $<3\times 10^{-6}$&0.000483\\
\hline
2 & 13 & 0.000389 & 162 & 0.000387 & $<3\times 10^{-6}$&0.000384\\
\hline
2 & 14 & 0.000315 & 186 & 0.000314 & $<4\times 10^{-6}$&0.000310\\
\hline
2 & 15 & 0.000259 & 222 & 0.000258 & $<5\times 10^{-6}$&0.000253\\
\hline
2 & 16 & 0.000215 & 258 & 0.000214 & $<5\times 10^{-6}$&0.000209\\
\hline
2 & 17 & 0.000181 & 278 & 0.000180 & $<6\times 10^{-6}$&0.000174\\
\hline
2 & 18 & 0.000153 & 314 & 0.000152 & $<7\times 10^{-6}$&0.000145\\
\hline
2 & 19 & 0.000131 & 345 & 0.000130 &$<8\times 10^{-6}$&0.000122\\
\hline
2 & 20 & 0.000113 & 402 & 0.000112 &$<9\times 10^{-6}$&0.000103\\
\hline
2 & 21 & 0.000098 & 435 & 0.000097 &$<11\times 10^{-6}$&0.000086\\
\hline
2 & 22 & 0.000086 & 461 & 0.000085 &$<12\times 10^{-6}$& 0.000073\\
\hline
2 & 23 & 0.000076 & 497 &0.000073 & $<12\times 10^{-6}$& 0.000061\\
\hline
2 & 24 & 0.000067 & 573 & 0.000065 & $<14\times 10^{-6}$& 0.000051\\
\hline
2 & 25 & 0.000059 & 597& 0.000058 & $<17\times 10^{-6}$ & 0.000041\\
\hline
2 & 26 & 0.000053 & 736&0.000051 & $<18\times 10^{-6}$&0.000033\\
\hline
2 & 27 & 0.000048 & 804 & 0.000046 & $<20\times 10^{-6}$ & 0.000026\\
\hline
3 & 1 & 0.020834 & 381& 0.020795 & $<3\times 10^{-6}$ & 0.020792\\
\hline
3 & 2 & 0.001083 & 423 &0.001077 &$<2\times 10^{-6}$&0.001075\\ 
\hline\hline

\end{tabular}

\end{table}

\section{Extending results to lattice polytopes}

First we have to express analytically the intuitively obvious fact
that if two polytopes cover each other (up to a small homothety) then
their coefficients of asymmetry cannot be far apart.

\begin{lemma}\label{lm:\pl:ca} Let $P'\subset P$ be two polytopes such that $P$ can be
covered by a translate of $\lambda P'$. Then, for any $\B w\in\inter
(P')$,
 \begin{equation}\label{eq:\pl:ca}\ca(\B w,P)\le |\lambda|\, \ca(\B w,P')+|\lambda|-1.\end{equation} \end{lemma}
 \smallskip{\it Proof.}  Assume that $|\lambda|>1$, for otherwise $P'=P$ and we are
home. Also, the case of $1$-dimensional polytopes is trivial.

Let $\B w_1,\B w_2\in\partial P$ be two points with $\B w\in[\B w_1,\B
w_2]$ and $\ca(\B w,P)=|\B w_1-\B w|\,:\,|\B w-\B w_2|$. Let $\partial
P'\cap [\B w_i,\B w]=\{\B w_i'\}$, $i=1,2$, where $[\B x,\B y]$
denotes the straight line segment between $\B x$ and $\B y$. Clearly,
 \begin{equation}\label{eq:\pl:t1}|\B w_1-\B w_2|=(\ca(\B w,P)+1)\,|\B w-\B w_2|\ge (\ca(\B
w,P)+1)\,|\B w-\B w_2'|.\end{equation}

As $\lambda P'$ covers $\{\B w_1,\B w_2\}$, there are $\B u_1,\B
u_2\in P'$ with 
 \begin{equation}\label{eq:\pl:t2}\B w_1-\B w_2=|\lambda|\,(\B u_1-\B u_2).\end{equation}

We can assume that $\B u_1,\B u_2\in\partial P'$. If $\B u_1=\B w_1'$
and $\B u_2=\B w_2'$, then we let $\B v=\B u_2$. Otherwise let $\B v$
be the (well-defined) point of intersection of $L(\B u_1,\B w)$ and
$L(\B u_2,\B w_2')$, where $L(\B x,\B y)$ denotes the line though the
points $\B x$ and $\B y$. As $\B v\in L(\B u_2,\B w_2')$ lies outside
of $\inter(P')$, we have
 $$\ca(\B w,P')\ge \frac{|\B u_1-\B w|}{|\B w-\B v|}=\frac{|\B u_1-\B
u_2| \,-\, |\B w-\B w_2'|}{|\B w-\B w_2'|},$$
 which implies the required by~\eqref{eq:\pl:t1} and~\eqref{eq:\pl:t2}.\qed \medskip

\smallskip\noindent{\bf Remark.}  Note that the bound in~\eqref{eq:\pl:ca} is sharp, as is demonstrated
e.g.\ by $C_d\subset |\lambda| C_d$ and $\B w=c \B 1$ with $0<c\le
\frac12$, where $C_d\subset\I R^d$ is the $0/1$-cube.\smallskip

Now we are ready to prove our result on lattice polytopes.

\begin{theorem}\label{th:\pl:ca-polytope} Let $l\ge 1$ be an integer and let $P\subset\I R^d$
be a lattice polytope with $I_l(P)\not=\emptyset$. Then there is $\B
w\in I_l(P)$ with
 \begin{equation}\label{eq:\pl:ca-polytope-th} \ca(\B w,P)\le \frac{8d}{\delta(2d,\frac87\,
l)}-1= 8d\cdot (8l+7)^{2^{2d+1}}-1.\end{equation}\end{theorem}
  \smallskip{\it Proof.}  Let $S=\conv\{\B v_0,\dots \B v_d\}\subset P$ be a simplex of
the maximum volume; we may assume that each $\B v_i$ is a vertex of
$P$.

Choose $\B u\in I_l(P)$. Let $\B u_1\in\inter(S)$ be any vertex and
let $\B u_2$ be the point of intersection of the ray $\{\B u_1+\lambda
(\B u-\B u_1)\mid \lambda\ge 0\}$ with the boundary of $P$. The vertex
$\B u_2$ lies in the interior of some face which is spanned by at most
$d$ vertices of $P$. Hence $\B u\in\relinter([\B u_1,\B u_2])$ can be
represented as a positive convex combination of $n\le 2d+1$ vertices
of $P$ including all vertices of $S$, let us say $\B u=\sum_{i=0}^n
\alpha_i' \B v_i$ with $\sum_{i=0}^n \alpha_i'=1$ and each
$\alpha_i'>0$.

Choose a vertex $\B w\in I_l(P')$ and a representation $\B
w=\sum_{i=0}^n \alpha_i \B v_i$ with $\sum_{i=0}^n \alpha_i=1$
maximizing $\min_{0\le i\le n} \alpha_i$. Denote this maximum by $m(\B
w)>0$. The argument of Theorem~\ref{th:\pl:ca-simplex} shows that $m(\B w)\ge
\delta(n-1,\frac87\,l)/8\ge \delta(2d,\frac87\,l)/8$.

The polytope $P'$ can be represented as a projection of an $n$-simplex
$S_n$ such that $\B w$ is the image of $\B v\in\inter(S_n)$ with
$m_{S_n}(\B v)=m(\B w)$. Now, it is easy to see that a linear mapping
cannot increase the coefficient of asymmetry; hence 
 $$\ca(\B w,P')\le\ca(\B v,S_n)=\frac{1-m(\B w)}{m(\B w)}.$$

It is known that $P\subset (-d)S+(d+1)\B s$, where $\B s$ is the
centroid of $S$, see e.g.~\cite[Theorem~3]{lagarias+ziegler:91}. By
Lemma~\ref{lm:\pl:ca}, we obtain
 $$\ca(\B w,P)\le d\, \ca(\B w,P')+d-1 \le d\,\frac{1-m(\B w)}{m(\B
w)}+d-1=\frac{d}{m(\B w)}-1,$$
 which give the required by~\eqref{eq:\pl:LZ4}.\qed \medskip

\smallskip\noindent{\bf Remark.}  The bound~\eqref{eq:\pl:ca-polytope-th} is much worse
than~\eqref{eq:\pl:ca-simplex-i}; the reason is that we may have to approximate
$2d$-tuples of numbers in Lemma~\ref{lm:\pl:LZ91}. Unfortunately, we cannot
guarantee that $P'$ has much fewer than $2d+1$ vertices, as e.g.\
$P=\conv\{\pm\B e_1,\dots, \pm\B e_{d-1}, (k+1)l \B e_d\}$
demonstrates. Perhaps, one can show that any such example cannot be
extremal for our problem and thus improve on~\eqref{eq:\pl:ca-polytope-th}.\smallskip

\smallskip\noindent{\bf Remark.}  It should be possible to generalize Theorem~\ref{th:\pl:ca-simplex}
and~\ref{th:\pl:ca-polytope} by proving the existence of a number
$b=b(d,l,m)>0$ such that any lattice polytope $P$ contains $m$
distinct points in $I_l(P)$ (provided $|I_l(P)|\ge m$) with
coefficient of asymmetry of each being at least $b$. The idea of the
proof is the following. If $P$ is a simplex, take distinct $\B
w_1,\dots,\B w_m\in I_l(P)$ with with largest $m_P$'s. Now each jump
of $\B w_i$ either does not increase $m_P(\B w_i)$ or maps $\B w_i$
into some other $\B w_j$. We are done if we can show that if $m_P$ is
very small then there are at least $m$ distinct jumps increasing
it. The latter claim would be achieved by rewriting the proof of
Lemma~\ref{lm:\pl:LZ91}, so that in the conclusion we have at least $m$
suitable $(n+1)$-tuples of integers. To extend the claim to general
lattice polytopes, observe that $m$ vertices in $I_l(P)$ can be
represented each as a positive combination of $d+1$ vertices of a
max-volume symplex and at most $md$ other vertices of $P$ and follow
the argument of Theorem~\ref{th:\pl:ca-polytope}. We do not see any principal
difficulties arising here, but it would take too much space to write
the complete proof, so we restrict ourselves to this little
observation only.\smallskip

\section{Volume of lattice simplices}

For simplices we have a better method (than applying~\eqref{eq:\pl:M39}) for
bounding volume which appears in~\cite[Theorem~3.4]{hensley:83} (see
also~\cite[Lemma~2.3]{lagarias+ziegler:91}). Let us reproduce
this simple argument here.

\begin{lemma}\label{lm:\pl:H83} Let $S=\conv\{\B v_0,\dots,\B v_d\}$ be any simplex and let 
$\B w\in I_l(S)$ have barycentric coordinates $(\alpha_0,\dots,\alpha_d)$.
Then\begin{equation}\label{eq:\pl:H83}
 \vol(S)\le \frac{l^d}{d!\times \alpha_1\times\alpha_2\times\cdots\times
\alpha_d}\,|I_l(S)|.\end{equation} \end{lemma} 
 \smallskip{\it Proof.}  Consider the region
 $$X=\{\B w+{\textstyle \sum_{i=1}^d \beta_i(\B v_i-\B v_0)} \mid |\beta_i|\le
\alpha_i\ 1\le i\le d\}.$$
 It is a centrally symmetric parallelepiped around the vertex $\B w\in
l\I Z^d$ with volume $\vol(X)=d!\,\vol(S) \prod_{i=1}^d
(2\alpha_i)$. We have to show that the volume of $X$ cannot exceed
$(2l)^d\,|I_l(S)|$. If this is not true, then $X$ contains (besides
$\B w$) at least $|I_l(S)|$ pairs of vertices $\B w\pm\B u\in l\I Z^d$
by Corput's theorem~\cite{corput:35} and, clearly, at least one vertex
of each such pair lies within $I_l(S)$, which is a contradiction.\qed \medskip

Now we can deduce the following result.

\begin{theorem}\label{th:\pl:sdkl} For any $k\ge 1$, the inequality~\eqref{eq:\pl:sdkl-i} holds.\end{theorem}
 \smallskip{\it Proof.}  Let $S\subset \I R^n$ be lattice simplex with $|I_l(S)|=k$.
By Theorem~\ref{th:\pl:ca-simplex} there is
$\B w\in I_l(S)$ with $m_S(\B w)\ge \gamma=\smdl$. Assume
$\alpha_0\le\cdots\le \alpha_d$. Then it is easy to see that
$\prod_{i=1}^d \alpha_i\ge \gamma^{d-1}\, (1-\gamma d)$. The
claim now follows from~\eqref{eq:\pl:H83}.\qed \medskip

\section{$s(d,k,l)$ for small $d$ and $l$}

As we have already mentioned, $s(2,k,1)$ was computed by
Scott~\cite{scott:76}. The simplex $S_{2,k,l}$ shows that $s(2,k,l)\ge
l(l+1)^2 (k+1)/2$. Upper bounds on $s(2,k,l)$ can be obtained
by applying~\eqref{eq:\pl:H83} to the lower bounds on $\beta(2,l)$ from
Table~\ref{beta-cplex}. But even if we knew $\beta(2,l)$ exactly, the
best upper bound on $s(2,k,l)$ that this method would give is
$l^5k/2+O(l^4)k$, so there would still be an uncertainty
about $s(2,k,l)$.

Also, an interesting problem is the determination of $s(3,k,1)$. The
simplex $S_{3,k,1}$ shows that $s(3,k,1)\ge
6(k+1)$. Theorem~\ref{th:\pl:sdkl} gives, already for such small $d$, very
bad bounds. However, there is a very simple argument, following
the lines of Section~\ref{small-beta} and proving that
 \begin{equation}\label{eq:\pl:s3k1} s(3,k,1)\le \frac{29791}{2112}\cdot k< 14.106\cdot k.\end{equation}

Given a lattice simplex $S\subset\I R^3$, we can deduce as before that
the barycentric coordinates $(\alpha_0,\dots,\alpha_3)$ of a lattice
vertex maximizing $m_S$ satisfy $2\alpha_3-1\le \alpha_0$,
$3\alpha_2-1\le \alpha_0$ and either $4\alpha_2-1\le\alpha_0$ or
$4\alpha_3-2\le \alpha_0$. These inequalities do not guarantee yet
that $\alpha_0>0$, but they guarantee that $\alpha_1\ge \frac2{31}$
and, as it is routine to see, that $\alpha_1\alpha_2\alpha_3\ge
\frac2{31}\times\frac {11}{31}\times \frac{16}{31}$, which
implies~\eqref{eq:\pl:s3k1} by~\eqref{eq:\pl:H83}.

Of course, we could write more equations on the $\alpha$'s, but this
method would not lead to the best possible bound. For example, the
simplex $S_{3,1,1}$ shows that we cannot guarantee a vertex in
$I_1(S)$ with $\alpha_1\alpha_2\alpha_3 >
\frac12\times\frac13\times\frac1{12}$, so the best bound we would hope
to obtain this way is $s(3,k,1)\le 12k$ only.

\section{Lawrence's Finiteness Theorem}\label{lawrence}

A result of Lawrence~\cite[Lemma~5]{lawrence:91} implies that there is
$\gamma=\gamma(d)>0$ such that, for any lattice simplex $S=\conv\{\B
v_0,\dots,\B v_d\}$, the set $I_1(S)$ (if non-empty) contains a vertex
$\B w$ with $\alpha_0\ge\gamma$, where $(\alpha_0,\dots,\alpha_d)$ are
the barycentric coordinates of $\B w$. This directly follows from our
Theorem~\ref{th:\pl:ca-simplex} but unfortunately I could not find a simple
argument giving the converse implication. And, in fact, the
corresponding extremal functions are different: for example,
$\beta(2,1)\le 1/8$ while it is claimed in~\cite[p.~439]{lawrence:91}
that we can take $\gamma(2)=1/6$.

However, one can deduce from~\cite{lawrence:91} that $\beta(d,1)>0$
using the following simple modification of Lawrence's proof. For
$i\in\I N$ let $U_i=\{\B w\in\Delta_d\mid \ca(\B w,\Delta_d)<i\}$,
where $\Delta_d=\conv\{\B 0,\B e_1,\dots,\B e_d\}$. Clearly,
$\cup_{i=d+1}^\infty
U_i=\inter(\Delta_d)$. By~\cite[Theorem~3]{lawrence:91} there exists
$i$ such that for any $\B w\in\inter(\Delta_d)$ there are $j\in\I N$
and $\B u\in \I Z^d$ with with $j\B w+\B u\in U_i$. We claim that
$\beta(d,1)\ge \frac1{i+1}$. Indeed, let $S\subset \I R^d$ be any
lattice simplex and let $\B v\in I_1(S)$. Choose any affine function
$f:\I R^d\to\I R^d$ with $f(S)=\Delta_d$ and let $\B w=f(\B
v)\in\inter(\Delta_d)$. Given $\B w$, choose the corresponding $j\in\I
N$ and $\B u\in \I Z^d$. It is easy to check that $\B v'=f^{-1}(j\B
w+\B u)$ belongs to $I_1(S)$ and satisfies $m_S(\B v')\ge
\frac1{i+1}$.

But as it has already been remarked, Lawrence's argument does not give
any explicit bound.

\section*{Acknowledgments}

I am grateful to G\"unter Ziegler for drawing my attention to the
problem and for very useful discussions and suggestions at all
stages of this research. Also, I want to thank Marc Pfetsch and Julian
Pfleifle for help with programming and the whole "Discrete Geometry"
group of TU-Berlin for their hospitality.

\newpage
\pagestyle{empty}

\section*{C source code}

\small
\renewcommand{\baselinestretch}{0.84}
\hoffset= -1.5 true cm

\begin{verbatim}

/*

program for proving lower bounds on beta(d,l), the function defined in
e-print O.Pikhurko "Lattice Points inside Lattice Polytopes" at arXiv.org

Copyright (C) 2000 Oleg Pikhurko

This program is free software; you can redistribute it and/or
modify it under the terms of the GNU General Public License
as published by the Free Software Foundation.

This program is distributed in the hope that it will be useful,
but WITHOUT ANY WARRANTY; without even the implied warranty of
MERCHANTABILITY or FITNESS FOR A PARTICULAR PURPOSE. See the
GNU General Public License for more details:
http://www.gnu.org/copyleft/gpl.html


Purpose

For a lattice simplex $S\subset R^d$ and $x\in S$, let $m(x)$ be the
value of the smallest of the $d+1$ barycentric coordinates of
$x$. There is $\beta(d,l)>0$ such that if $S$ contains at least one
vertex of $lZ^d$ in its interior, then it contains such a vertex $x$
with $m(x)\ge \beta(d,l)$. The program computes lower bounds on
$\beta(d,l)$ for small $d$ and $l$.


Algorithm

Main idea: let $x\in S$ be a pair (almost) attaining $\beta(d,l)$. Let
$(a_0\le\dots\le a_d)$ be the barycentric coordinates of $x$. For any
numbers $p_0,\dots,p_d\in N$, there is $i\in[0,d]$ such that
$(lq+1)a_i-lp_i\le a_0$, where $q=\sum_{i=0}^d p_i$. (For otherwise we
have a contradiction, the vertex $x'=(q+1)x-\sum_{i=1}^d p_i v_i\in
lZ^d$ belongs to the interior of $S$ and has larger $m$, where $v_i$'s
are the vertices of $S$.)

That is, the $a_i$'s (for extremal $x$) must satisfy certain
linear inequalities. If these imply a lower bound on $a_0$ it is also
a lower bound on $\beta(d,l)$.

In the algorithm, we define our initial LP to consists of $0\le
a_0\le\dots\le a_d$ and $\sum_{i=0}^d a_i=1$. Then we repeat the
following. Given LP, let $(a_i)$ be a solution minimizing $a_0$. Try
to find $p_i$ with $(lq+1)a_i-lp_i>a_0$ which show that $(a_i)$ cannot
be the coordinates of extremal $x$, add the corresponding constraints
to our LP, and repeat.

Of all $p_i$'s, we choose the smallest in the lex order with minimum
$q$, which probably speeds convergence. Any better ideas?
 
As each time we have or-connected inequalities, we introduce binary
variables; of course, if $p_i=p_{i-1}$ then there is not need to
include the $i$th inequality. Also the $0$th equality is never
included because $(a_i)$ minimizes $a_0$ given LP, hence if
$(lq+1)a_0-lp_0>a_0$ for this LP, it will be true for any larger LP.


Usage

Prepare file "d-l-in.lp(s)" containing the lp program at which you
wish to start. (If the file does not exist, the program creates the
default initial lp.) Tip: rename file "d-l-last.lp(s)" as
"d-l-in.lp(s)" to start at the place, where the program finished last
time.

Run the program. Enter d l when asked. As each subsequent prompt,
enter 0 for non-iteractive mode, -1 to finish calculations, -2 to save
lp to the file "d-l-out.lp(s)" or a positive number which gives the
number of iterations to do till next prompt. At each termination of
the program (except by ^C), the last LP problem is saved to file
"d-l-last.lp(s)".

The program can be linked with either cplex (www.cplex.com) or lp_solve 
(ftp://ftp.ics.ele.tue.nl/pub/lp_solve/). With this problem lp_solve is 
unfortunately not very stable. Required libraries: 
  cplex, socket, m, nsl (if compiled with -DCPLEX);
  liblpk.a, m, l (if compiled with the -DLPS).


Bugs

If we compile with -DLPS and read the initial lp from file, then
change_lp causes segmentation fault. However, everything works fine if 
we compile with -DCPLEX or the initial problem is created by make_lp().

*/


#include <stdio.h>
#include <math.h>
#include <string.h>

/* library specific declarations; some of our function also depend
   on the library; these are placed at the end of the file */
#ifdef CPLEX

 #include <cplex.h>
 #define REAL double
 REAL tolerance=(REAL)1e-9; /* min tolerance allowed by CPLEX */
 REAL delta=(REAL)1e-9; /* initial upper bound on possible error on a[0]; 
                    for CPLEX it will be recomputed with each iteration */
 
 /* extension of all lp files */
 #define EXT ".lp"
 CPXENVptr env = NULL;
 CPXLPptr lp = NULL;
 int status = 0;

#elif defined LPS

 #include <lpkit.h>

 /* set different extension to avoid interference */
 #define EXT ".lps"
 lprec *lp;
 REAL delta=(REAL)1e-5; /* this is our guess for lp_solve */

#else

 #error Please compile with either -DCPLEX or -DLPS option

#endif




#define DMAX 9  /* size allocated for a[]; 
                   d must not exceed DMAX */
#define QMAX 1000 /* if next_p we get as far as q>=QMAX,
                     we assume that no approximation exists */
#define IMAX 600 /* the maximal number of iterations */

/* most important functions */
void get_a(void), read_a(void); /* get/(read from stdin) next solution a[] */
void next_p(void); /* compute approximation p[] */
void change_lp(void); /* add the corresponding constraints to lp */

/* auxiliary functions */
void print_a(REAL*), print_p(int*), print_all(void); /* various printing */
void die(char*); /* terminate the program */
void close_all(void); /* close all files, etc */
void init_lp(void); /* initialization */
void add_dl(char*); 
void find_jumps(void);
void p_string(char*,int);
void save_lp(char*);
void check_solution(void);
void change_lp_plain(void);
void change_lp_tricky(void);

int d=2,l=1;
REAL best=(REAL) 0.0; /* the best lower bound we can guarantee so far */
REAL e=(REAL) 1e-6; /* small constant to be used in crucial roundings,
                    we search for p[] with (l*q+1)a[i]-l*p[i] > a[0]+e  */
REAL a[DMAX+1], b[DMAX+1]; /* b[] = (l*q+1)a[]-l*p[] */
int p[DMAX+1]; /* a, b & p have d+1 elements each */
int q;  /* q=p[0]+...+p[d] */
int jp, jp_ind[DMAX+1]; /* the number and array of jumps in p[] */
int it; /* the number of iterations */
int  d_stop; /* gives number of iteration before next prompt; 
                if <=0 used for various actions */



int main()
{

  int next_stop=0;

  printf("Please enter d l ");
    scanf("%d %d",&d,&l);

  if (d>DMAX || d<=0 || l<=0)
    die("Wrong input values for d and l.");

  init_lp();

  for(it=0;it<IMAX;it++) {

    get_a(); /* read_a(); */

    check_solution(); /* can be commented out */

    if(best<a[0]-delta) /* we can improve our lower bound! */
      best=a[0]-delta;

    next_p();

    print_all();

    change_lp();

    if(it>=next_stop) {

      printf("\n\nEnter 0 (non-stop), -1 (end) -2 (save lp)\n");
      printf("or number (>0) runs till stop ");
      scanf("%d",&d_stop);

      if(d_stop==-1)
        die("Exiting...");

      if(d_stop==-2) {
        char lp_name[]="x-x-out" EXT;
        save_lp(lp_name);
      }

      if(d_stop==0) 
        next_stop=IMAX+1;
     
      if(d_stop>0)
        next_stop=it+d_stop;

    }

  }

  die("The maximum number of iterations is exceeded.");

  return 0;

}



/* reads a[] from stdout; for testing purposes */
void read_a(void)
{

  int i;
  REAL suma;

  printf("Enter a\n");

  for(suma=0,i=0;i<d;i++) {
    scanf("%lf", &a[i]);
    suma+=a[i];
  }

  a[d]=1-suma;

  if (a[0]<0 || a[d]<0)
    die("a[] is not feasible.");

}



/* for given a[] compute a smallest p[] such that min((lq+1)a[]-l*p[])>a[0].
   return 0 if no reasonably small P exists */
void next_p(void)
{

  int found=0,i,sump,diff;

  /*  printf("next_p: obtained the following vector:\n");
  print_a(a); */

  q=0; /* we try to find p[] with smallest q, so we interate on q */
  
  do {

    q++;


    sump=0;
    for(i=0;i<=d;i++) {

      /* the case i==0 receives a special treatment */
      p[i]=ceil(q*(a[i]-(i?0:delta)) + (a[i]-a[0])/l - 1 - e);

      if (p[i]<0) /* this may happen if a[i] is small */
        p[i]=0;
      sump+=p[i];
    }

    if (sump>=q) {

      found=1;

      if (sump>q) {  /* here different p[] are possible;
                        we take the lex-smallest p[] */ 
        diff=sump-q;
        for(i=0; i<=d && diff>0; i++)
          if(p[i]<=diff) {
            diff-=p[i];
            p[i]=0;
          } else {
            p[i]-=diff;
            diff=0;
          }
      }

    }

  } while (! found && q< QMAX);

 
  if (found){

    /* let us run small check; just in case */
    for(sump=0,i=0;i<=d;i++)
      sump+=p[i];

    if (sump != q)
      die("next_p: sump does not equal q");

    for(i=0;i<=d;i++) {
      b[i]=(l*q+1)*a[i]-l*p[i];
      if( p[i]<0 || ( b[i]<= a[0] && a[i]>0 ) )
        die("next_p: the sequence p[] is bad");
    }

  } else
    die("Suitable p cannot be found. Success!?");

}


/* creates a string "r_p[0]_..._p[d]_i" */  
void p_string(char s[], int i)
{

  int j;
  char temp_name[255];

    strcpy(s,"r");
    for(j=0;j<=d;j++) {
      sprintf(temp_name,"_%d",p[j]);
      strcat(s,temp_name);
    }
    sprintf(temp_name,"_%d",i);
    strcat(s,temp_name);

}



/* modifies our lp */
void change_lp(void)
{

  find_jumps(); /* compute the jump in p[] */

  if(!jp) /* a check just in case */
    die("No jumps in p[], which cannot happen at all.");
      
  change_lp_tricky();

  /* change_lp_plain(); */

  /* change_lp_plain() produces lp which are longer to solve, but
     whose binary variables are easier to understand (by humans);
     we keep the definition of change_lp_plain() for debugging */

}



/* compute the number jp and array jp_ind[] of jumps in p[]  */
void find_jumps(void)
{
  
  int i,p_prev;

  jp=0;
  p_prev=p[0];

  for(i=1;i<=d;i++)
    if (p[i]>p_prev) {
      jp_ind[jp++]=i;
      p_prev=p[i];
    }

}



/* prints d+1 array of reals */
void print_a(REAL s[])
{
 
  int i;

  for(i=0;i<=d;i++)
    printf("%19.17f ",s[i]);
  printf("\n");

}




/* prints d+1 array of integers */
void print_p(int s[])
{
 
  int i;

  for(i=0;i<=d;i++)
    printf("%d ",s[i]);
  printf("\n");

}



/* prints some global information */
void print_all()
{

 

  printf("\nGlobal variables:\n");
  printf("d=%d, l=%d, iteration=%d, q=%d, delta=%.2e\n",
         d,l,it,q,delta);

#ifdef CPLEX
  if(status) { /* get the error message corresponding to status */
      char  errmsg[1024];
      CPXgeterrorstring (env, status, errmsg);
      printf ("%s\n", errmsg);
  }
#endif

  printf("Best guaranteed lower bound so far = %19.17f\n", best);
  printf("Vectors a[] p[] & b[]:\n");
  print_a(a);
  print_p(p);
  print_a(b);

}



/* stop execution */
void die(char s[])
{

  printf("\n%s\n",s);
  print_all();

  close_all();
 
  exit(1);

}



/* replaces the 1st and 3d latters in a filename by the values of d and l */

void add_dl(char s[])
{
 
  s[0]='0'+d;
  s[2]='0'+l;

}


/* check if new a[] is compatible with p[]; lp_solve sometimes returns
   an unfeasible solution and check_solution fails - any ideas? */
void check_solution(void)
{

  int i,found=0;
  
  for(i=0;i<=d;i++)
    if((b[i]=(l*q+1)*a[i]-l*p[i]) <= a[0]+e)
      found=1;

  if (!found)
    die("check_solution failed.");

}





/* library specific definitions of functions: change_lp_plain(), 
                change_lp_tricky(), close_all(), save_lp() */




#ifdef CPLEX



/* Initialization */
void init_lp(void) 
{

  int i;

  char in_lp[]="x-x-in" EXT;
  FILE *lp_file;

  add_dl(in_lp);
  
  /* Initialize the CPLEX environment */
  env=CPXopenCPLEXdevelop(&status);
  if (env == NULL)
    die("CPXopenCPLEXdevelop failed.");

  /* sets extra space for columns */
  status=CPXsetintparam(env,CPX_PARAM_COLGROWTH,DMAX*IMAX);
  if(status)
    die("Cannot set CPX_PARAM_COLGROWTH.");

  /* sets extra space for rows */
  status=CPXsetintparam(env,CPX_PARAM_ROWGROWTH,(DMAX+1)*IMAX);
  if(status)
    die("Cannot set CPX_PARAM_ROWGROWTH.");

  /* sets extra space for non-zero elements;
     do we need this, given two settings above? */
  status=CPXsetintparam(env,CPX_PARAM_NZGROWTH,(DMAX+1)*(DMAX+1)*IMAX);
  if(status)
    die("Cannot set CPX_PARAM_NZGROWTH.");

  /* sets integrality tolerance to smallest allowable */
  status=CPXsetdblparam(env,CPX_PARAM_EPINT,tolerance);
  if(status)
    die("Cannot set CPX_PARAM_EPINT.");

  /* sets optimallity tolerance to smallest allowable */
  status=CPXsetdblparam(env,CPX_PARAM_EPOPT,tolerance);
  if(status)
    die("Cannot set CPX_PARAM_EPOPT.");

  /* sets feasibility tolerance to smallest allowable */
  status=CPXsetdblparam(env,CPX_PARAM_EPRHS,tolerance);
  if(status)
    die("Cannot set CPX_PARAM_EPRHS.");
 
  /* create lp environment */
  lp=CPXcreateprob(env,&status,"plattice");
  if (lp == NULL)
    die("CPXcreateprob failed.");

  if((lp_file=fopen(in_lp, "r"))==NULL) {
  
    REAL obj[DMAX+1]; /* objective function */
    char colname[DMAX+1][4]; /* names of variables, "a0", "a1", etc */
    char* colname_p[DMAX+1];
    REAL rhs[DMAX+1],rmatval[3*DMAX+1];
    int rmatbeg[DMAX+2],rmatind[3*DMAX+1];
    char sense[DMAX+1];

    printf("File %s does not exists.\n",in_lp);
    printf("Creating the basic initial lp.\n");

    for(i=0;i<=d;i++) {
      obj[i]=(REAL)0.0;
      colname_p[i]=strcpy(colname[i],"a0");
      colname[i][1]+=i;
    }

    obj[0]=(REAL)1.0;

    /* add real, >=0 (by default) variables a0,a1,... */
    status=CPXnewcols(env,lp,d+1,obj,NULL,NULL,NULL,colname_p);
    if (status)
      die("CPXnewcols failed.");


    for(i=0;i<d;i++) { /* inequality -a[i] + a[i+1] >= 0 */
      rhs[i]=(REAL)0.0;
      sense[i]='G';
      rmatbeg[i]=2*i;
      rmatind[2*i]=i;
      rmatval[2*i]=(REAL)-1.0;
      rmatind[2*i+1]=i+1;
      rmatval[2*i+1]=(REAL)1.0;
    }

    /* equation a0+a1+...=1 */
    rhs[d]=(REAL)1.0;
    sense[d]='E';
    rmatbeg[d]=2*d;
    for(i=0;i<=d;i++) {
      rmatind[2*d+i]=i;
      rmatval[2*d+i]=(REAL)1.0;
    }
      
    /* add these basic constraints */
    status=CPXaddrows(env,lp,0,d+1,3*d+1,rhs,sense,
                      rmatbeg,rmatind,rmatval,NULL,NULL);
    if (status)
      die("CPXaddrows failed.");
    
    } else {

    fclose(lp_file);

    status=CPXreadcopyprob(env,lp,in_lp, NULL);
    if (status)
      die("CPXreadcopyprob failed.");

  }

}



/* compute a[] as the next solution; */
void get_a(void)
{


  /* solve the lp */

  status=CPXmipopt(env, lp);
  if (!status) {

    status = CPXgetmipx(env, lp, a, 0, d);
    if (status) 
      die("CPXgetmipx failed\n");

  } else { /* CPXmipopt failed; maybe there are no bin variables yet */

    /*    printf("%d\n",status); */

    status=CPXprimopt(env,lp);
    if(status)
      die("Both CPXmipopt and CPXprimopt failed.");
    status=CPXgetx(env,lp,a,0,d);
    if(status)
      die("CPXgetx failed.");

  }

}



/* add constraints corresponding to p[] to the lp;
   a straightforward version */
void change_lp_plain(void)
{

  int j,i;
  int cur_numcols; /* the number of columns (variables) in lp */
  /* variables for CPXaddcols */

  REAL ub[DMAX];

  /* variables for CPXaddrows; 
     we may need as much as d+1 equations */
  REAL rhs[DMAX+1];
  int rmatbeg[DMAX+2];
  int rmatind[4*DMAX]; 
  char sense[DMAX]; /* types of constraints */
  char row_name[DMAX+1][255]; 
  char* row_pointer[DMAX]; /* symbolic names of rows */ 
  REAL rmatval[4*DMAX]; /* our equations have <=4*d non-zero coefficients */
  char ctype[DMAX]; /* types of variables */

  /* first we add jp binary variables z[0],...,z[jp-1] */
  /* for j=1 we add an unecessary binary variable,
     but it is fixed to 1 anyway and brings no harm;
     we must also specify ub[j]=1 because otherwise
     CPXmipopt returns 3002. */
  for(j=0;j<jp;j++) {
    ctype[j]='B';
    ub[j]=(REAL)1.0;
  }
  status=CPXnewcols(env,lp,jp,NULL,NULL,ub,ctype,NULL);
  if(status)
    die("CPXnewcols failed.");

  /* get the number of columns (variables) in lp */
  cur_numcols = CPXgetnumcols (env, lp);
  if(cur_numcols<=0)
    die("CPXgetnumcols failed.");

  /* now out added variables have numbers
     from cur_numcols-jp to cur_numcols-1 */

  /* now we add jp constraints of the form 
     -a[0]+(l*q+1)*a[i]+(l*q+1) z[j]<= l*p[i]+(l*q+1), where i=jp_ind[j] */

  if(delta<tolerance*(3+2*l*q)) /* recompute the max abs mistake */
    delta=tolerance*(3+2*l*q);

  for(j=0;j<jp;j++) {

    i=jp_ind[j];

    rhs[j]=(REAL) l*p[i]+(l*q+1);
    rmatbeg[j]=3*j; /* three non-zero coeff. per row */
    rmatval[3*j]=(REAL) -1; 
    rmatind[3*j]=0; /* x[0] */
    rmatval[3*j+1]=(REAL) l*q+1;
    rmatind[3*j+1]=i; /* x[i] */
    rmatval[3*j+2]=(REAL) l*q+1; 
    rmatind[3*j+2]=cur_numcols-jp+j; /* z[j] */
    sense[j]='L';
    p_string(row_name[j],i);
    row_pointer[j]=row_name[j];

  }

  /* plus the constraint z[0]+...+z[jp-1]>=1; */

  if(delta<tolerance*jp) /* recompute the max abs mistake */
    delta=tolerance*jp;

  rhs[jp]=(REAL)1.0;
  sense[jp]='G';
  rmatbeg[jp]=3*jp;
  rmatbeg[jp+1]=4*jp;
  p_string(row_name[jp],0);
  row_pointer[j]=row_name[j];

  for(j=0;j<jp;j++) {
    rmatval[3*jp+j]=(REAL) 1;
    rmatind[3*jp+j]=cur_numcols-jp+j;
  }

  status=CPXaddrows(env,lp,0,jp+1,4*jp,rhs,sense,
                    rmatbeg,rmatind,rmatval,NULL,row_pointer);
  if(status)
    die("CPXaddrows failed.");

}





/* add constraints corresponding to p[] to the lp; uses the minimum
   nnumber of binary variables; a bit tricky */
void change_lp_tricky(void)
{

  /* variables for CPXaddrows; 
     we may need as much as d equations */
  REAL rhs[DMAX];
  int rmatbeg[DMAX+1];
  int rmatind[DMAX*(DMAX+2)]; /* non-zero elements in added equations */
  char sense[DMAX]; /* types of constraints */
  char row_name[DMAX][255]; 
  char* row_pointer[DMAX]; /* symbolic names of rows */ 
  REAL rmatval[DMAX*(DMAX+2)];
 
  int j, jp_log=0, j2=1;
  int cur_numcols; /* the number of columns (variables) in lp */

  /* we add constraints saying that there is j<jp such that
     (l*q+1)a[i]- a[0]<= l*p[i], where i=jp_ind[j];
     to do so we introduce jp_log binary variables, where 
     jp_log is the smallest number with $2^jp_log>= jp.

     For each of 2^jp_log possible combinations of binaries all
     inequalities are vacuous except one, which is precisely one of
     our desired inequalities. We do not explain how such a system
     is constructed but encourage the reader to check the algorithm
     and to run the program and look at the resulting lp */

  while(j2 < jp) {
    j2*=2;
    jp_log++;
  }

  if(delta<tolerance*(l*q+1)*(jp_log+1))  /* recompute the max mistake */
    delta=tolerance*(l*q+1)*(jp_log+1);

  if (jp_log>0) {

    /* variables for CPXaddcols */
    REAL ub[DMAX];
    char ctype[DMAX]; /* types of variables */

    /* first we add jp_log binary variables
       we must also specify ub[j]=1 because otherwise
       CPXmipopt returns 3002. */

    for(j=0;j<jp_log;j++) {
      ctype[j]='B';
      ub[j]=(REAL)1.0;
    }

    status=CPXnewcols(env,lp,jp_log,NULL,NULL,ub,ctype,NULL);
    if(status)
      die("CPXnewcols failed.");

  }

  /* now we add constraints, the case jp_log==0 inclusive */

  /* get the number of columns (variables) in lp */
  cur_numcols = CPXgetnumcols (env, lp);
  if(cur_numcols<=0)
    die("CPXgetnumcols failed.");
  /* our added variables have numbers
     from cur_numcols-jp_log to cur_numcols-1 */

  for(j=0;j<jp;j++) {

    int jj,jj2,i;

    i=jp_ind[j];

    sense[j]='L';

    rmatbeg[j]=j*(2+jp_log); /* 2+jp_log non-zero entries per row */

    rmatind[j*(2+jp_log)]=0; /* coefficient at a0 */
    rmatval[j*(2+jp_log)]=(REAL)-1.0;

    rmatind[j*(2+jp_log)+1]=i; /* coefficient at a[i] */
    rmatval[j*(2+jp_log)+1]=(REAL) (l*q+1);

    p_string(row_name[j],i);
    row_pointer[j]=row_name[j];

    rhs[j]=(REAL) l*p[i];

    for(jj=0,jj2=1; jj<jp_log; jj++,jj2*=2) {


      /* coefficient at jj'th binary variable */ 
      rmatind[j*(2+jp_log)+2+jj]=cur_numcols-jp_log+jj;

     
      if( ( (unsigned int) 1 << jj) & j ) {
        
        rhs[j]+=(REAL)l*q+1;
        rmatval[j*(2+jp_log)+2+jj]=(REAL)+l*q+1;
      } else {
        if (j+jj2>=jp)
          rmatval[j*(2+jp_log)+2+jj]=(REAL)0.0;
        else
          rmatval[j*(2+jp_log)+2+jj]=(REAL)-l*q-1;
      }
    }

  }

  status=CPXaddrows(env,lp,0,jp,jp*(2+jp_log),rhs,sense,
                    rmatbeg,rmatind,rmatval,NULL,row_pointer);
  if(status)
    die("CPXaddrows failed.");

}



/* close log-files, cplex and write the current program to file */

void close_all(void)
{

  char out_lp[]="x-x-last" EXT;

  save_lp(out_lp);
  CPXwriteprob(env,lp,out_lp, NULL);
  CPXfreeprob(env, &lp);
  CPXcloseCPLEX(&env);

}



/* saves the lp to file; we cannot here call die()
   because save_lp may be called from it */

void save_lp(char s[])
{

  add_dl(s);
  printf("Writing current LP to file %s ... ",s);
  status=CPXwriteprob(env,lp,s, NULL);

  if(status)
    printf("failed, status=%d.\n", status);
  else
    printf("ok.\n");

}




#elif defined LPS


/* read/create the initial lp */
void init_lp(void)
{

  nstring col_name;
  char in_lp[]="x-x-in" EXT;
  FILE *lp_file;
  
  add_dl(in_lp);

  if((lp_file=fopen(in_lp, "r"))==NULL) {

    int i,j;
    REAL row[DMAX+2];

    printf("File %s does not exists.\n",in_lp);
    printf("Creating the basic initial lp.\n");

    lp=make_lp(0,d+1);
    if(lp==NULL)
      die("make_lp failed.");
    
    for(i=0;i<=d;i++) { /* give variables names a0,a2,... */
      strcpy(col_name,"a0");
      col_name[1]+=i;
      set_col_name(lp,i+1,col_name);
    }

    /* constraint -a[i-1]+a[i]>=0 */
    /* add_constraint starts indexing with 1 and sometimes
       the 0th entry is also used for other purposes; so we
       have to modify all indexing accordingly */

    for(i=1;i<=d;i++) {
      for(j=0;j<=d+1;j++)
        row[j]=(REAL) 0.0;
      row[i]=(REAL) -1.0;
      row[i+1]=(REAL) 1.0;
      add_constraint(lp,row,GE,(REAL) 0.0);
    }

    /* a[0]+...+a[d]=1 */
    for(j=1;j<=d+1;j++)
      row[j]=(REAL) 1.0;
    add_constraint(lp,row,EQ, (REAL) 1.0);
      
    /* set objective function */
    
    row[1]=(REAL) 1.0;
    for(j=2;j<=d+1;j++)
      row[j]=(REAL) 0.0;
    
    set_obj_fn(lp,row);

    print_lp(lp);

  } else {

    /* NB: keep lp_file open */
    lp=read_lp_file(lp_file,FALSE,"plattice");
    fclose(lp_file);

    if(lp==NULL)
      die("read_lp_file failed.");

  }

}


/* compute a[] as the first d+1 variables of the solution; */
void get_a(void)
{

  int i;

  if(solve(lp))
    die("solve failed.");

  for(i=0;i<=d;i++)    
    a[i]=lp->best_solution[lp->rows+1+i];

}


/* add constraints corresponding to p[] to the lp;
   a straightforward version */
void change_lp_plain(void)
{

  nstring row_name; /* ! lp_solve allows at most 24 characters */

  int j,i,jj;

  REAL col[(DMAX+1)*(IMAX+1)];
  /* (DMAX+1)*(IMAX+1) is max possible number of rows/cols in lp */

  for(j=0;j<=lp->rows;j++)
    col[j]=(REAL) 0.0;

  /* first we add jp binary variables z[0],...,z[jp-1] */
  /* for j=1 we add an unecessary binary variable,
     but it is fixed to 1 anyway and brings no harm; */

  for(j=0;j<jp;j++) {
    add_column(lp,col);
    set_int(lp,lp->columns,TRUE);  /* the variable is integer */
    set_upbo(lp,lp->columns,(REAL) 1.0);  
    set_lowbo(lp,lp->columns,(REAL) 0.0);  
  }

  /* now we add jp constraints of the form 
     -a[0]+(l*q+1)a[i]+(l*q+1) z[j]<= l*p[i]+(l*q+1), where i=jp_ind[j] */
  for(j=0;j<jp;j++) {
    i=jp_ind[j];
    for(jj=0;jj<=lp->columns;jj++)
      col[jj]=(REAL) 0.0;
    col[1]=(REAL) -1.0; /* at a[0] */
    col[i+1]=(REAL) l*q+1; /* at a[i] */
    col[lp->columns-jp+j+1]= (REAL) l*q+1; /* at z[j] */
    add_constraint(lp,col,LE,(REAL) l*p[i]+(l*q+1));
    p_string(row_name,i);
    set_row_name(lp,lp->rows,row_name);
  }

  /* plus the constraint z[0]+...+z[jp-1]>=1 */
  for(jj=0;jj<=lp->columns;jj++)
    col[jj]=(REAL) 0.0;
  for(j=1;j<=jp;j++)
    col[lp->columns-jp+j]=(REAL) 1.0;
  add_constraint(lp,col,GE,(REAL) 1.0);
  p_string(row_name,0);
  set_row_name(lp,lp->rows,row_name);
      
}


/* add constraints corresponding to p[] to the lp; uses the minimum
   nnumber of binary variables; a bit tricky */
void change_lp_tricky(void)
{

  int j, jp_log=0, j2=1;

  nstring row_name; /* ! lp_solve allows at most 24 characters */
 
  REAL col[(DMAX+1)*(IMAX+1)];
  /* (DMAX+1)*(IMAX+1) is max possible number of rows/cols in lp */

  for(j=0;j<=lp->rows;j++)
    col[j]=(REAL) 0.0;

  /* we add constraints saying that there is j<jp such that
     (l*q+1)a[i]- a[0]<= l*p[i], where i=jp_ind[j];
     to do so we introduce jp_log binary variables, where 
     jp_log is the smallest number with $2^jp_log>= jp.

     For each of 2^jp_log possible combinations of binaries all
     inequalities are vacuous except one, which is precisely one of
     our desired inequalities. We do not explain how such a system
     is constructed but encourage the reader to check the algorithm
     and to run the program and look at the resulting lp */

  while(j2 < jp) {
    j2*=2;
    jp_log++;
  }
  
  /* first we add jp_log binary variables */
  for(j=0;j<jp_log;j++) {
    add_column(lp,col);
    set_int(lp,lp->columns,TRUE);  /* the variable is integer */
    set_upbo(lp,lp->columns,(REAL) 1.0);  
    set_lowbo(lp,lp->columns,(REAL) 0.0);  
  
  }

  /* our added variables have numbers
     from lp->columns-jp_log to lp->columns-1 */

  /* add constraints now */

  for(j=0;j<jp;j++) {

    int jj,jj2,i,jj_beg;
    REAL rhs;

    for(jj=0;jj<=lp->columns;jj++)
      col[jj]=(REAL) 0.0;

    i=jp_ind[j];

    col[1]=(REAL) -1.0; /* coefficient at a0 */
    col[i+1]=(REAL) l*q+1; /* coefficient at a[i] */

    rhs=(REAL) l*p[i]; /* rhs as it were without binary variables */

    for(jj=0,jj2=1; jj<jp_log; jj++,jj2*=2) {

      /* coefficient at jj'th binary variable */ 
      jj_beg=lp->columns-jp_log+jj+1;
      if( ( (unsigned int) 1 << jj) & j ) {     
        rhs+=(REAL)l*q+1;
        col[jj_beg]=(REAL) l*q+1;
      } else {
        if (j+jj2>=jp)
          col[jj_beg]=(REAL) 0.0;
        else
          col[jj_beg]=(REAL) -l*q-1;
      }
    }

    add_constraint(lp, col, LE, rhs);
    p_string(row_name,j);
    set_row_name(lp,lp->rows,row_name);
  }

}


/* write out d-l-last.lp(s) */ 
void close_all(void)
{

  char out_lp[]="x-x-last" EXT;

  save_lp(out_lp);

} 



/* saves the lp to file; we cannot here call die()
   because save_lp may be called from it */

void save_lp(char s[])
{

  FILE *lp_file;

  add_dl(s);
  printf("Writing current LP to file %s ... ",s);
 
  if((lp_file=fopen(s,"w"))==NULL)

    printf("failed.\n");

  else {

    write_LP(lp,lp_file); /* purify complains here about "uninitialized memory 
                             read", though the code seems clean to me */

    if(fclose(lp_file))
      printf("failed.\n");
    else
      printf("ok\n");

  }

}


#endif

\end{verbatim}

\end{document}